\documentclass[reqno,12pt]{amsart}
\usepackage{latexsym, amsfonts,mathrsfs, amsmath, amssymb, amscd, epsfig}
\usepackage{mathrsfs}

\newtheorem{thm}{Theorem}

\newtheorem{prop}[thm]{Proposition}
\newtheorem{lem}[thm]{Lemma}
\newtheorem{cor}[thm]{Corollary}
\newtheorem{rmk}{Remark}

\newenvironment{pf}{{\noindent \it \bf Proof:}}{{\hfill$\Box$}\\}

\def\Div{\text{div}}
\def\mcB{\mathcal{B}}
\def\mcD{\mathcal{D}}
\def\mcE{\mathcal{E}}
\def\mcF{\mathcal{F}}
\def\mcM{\mathcal{M}}
\def\msM{\mathscr{M}}
\def\mbR{\mathbb{R}}
\def\mbZ{\mathbb{Z}}
\def\mcS{\mathcal{S}}

\def\bB{\bar{B}}
\def\ckB{\check{B}}
\def\bu{\bar{u}}
\def\bv{\bar{v}}
\def\bm{\bar{m}}
\def\brho{\bar{\rho}}

\def\tOmegaN {\tilde{\Omega}_N}

\def\tH{\tilde{H}}
\def\tc{\tilde{c}}
\def\tB{\tilde{B}}
\def\tm{\tilde{m}}
\def\trho{\tilde{\rho}}
\def\hH{\hat{H}}
\def\hPsi{\hat{\Psi}}
\def\bpsi{\bar{\psi}}
\def\hc{\hat{c}}
\def\hm{\hat{m}}
\def\hu{\hat{u}}
\def\hv{\hat{v}}
\def\hrho{\hat{\rho}}

\begin{document}

\title[ ]{Existence of Steady Subsonic Euler Flows
 through Infinitely Long Periodic Nozzles}

\author[ ]{Chao Chen}
\address{The
Institute of Mathematical Sciences and department of mathematics, The Chinese University of Hong
Kong, Shatin, Hong Kong}
\email{cchen@math.cuhk.edu.hk}
\author[ ]{Chunjing Xie}
\address{Department of mathematics, University
of Michigan, 530 Church Street, Ann Arbor, MI 48109-1043 USA.}
\email{cjxie@umich.edu}


\begin{abstract} In this paper, we study the global existence of
steady subsonic Euler flows through infinitely long  nozzles which are periodic in $x_1$ direction with the period $L$. It is shown that when the
variation of Bernoulli function at some given section is
small and mass flux is in a suitable regime,
there exists a unique global subsonic flow
in the nozzle. Furthermore, the flow is also periodic in $x_1$ direction with the period $L$. If, in particular, the Bernoulli function is a constant, we also get the existence of subsonic-sonic flows when the mass flux takes the critical value.
\end{abstract}

\maketitle

\section{Introduction and Main Results}\label{Introduction}
The study on subsonic and transonic flows in nozzles has grown
enormously in recent years. Subsonic and subsonic-sonic potential
flows in infinitely long nozzles were studied in \cite{XX1,
XX2,LI, DuXinYan}.  For full compressible Euler equations, Xie and
Xin in \cite{XX3} showed the existence of global subsonic flow in
an infinitely long nozzle which tends  to be flat at far field.
The key point in \cite{XX3} is to transform the system of Euler
equations into a second order equation of stream function. The
careful energy estimates  give far field behavior and uniqueness
of flows. The idea in \cite{XX3} was generalized to subsonic Euler
flows in axially symmetric nozzles \cite{DD}. Subsonic and
subsonic-sonic potential flows past a body were studied in
\cite{Bers1,Bers2, CDSW} and references therein. Subsonic Euler
flows with nonzero vorticity in half space was investigated in
\cite{ChenJun}. Subsonic flows were also studied as a part of
stability of transonic shock problem, see \cite{CF1, CF2, CCF,
CCS, ChenS, ChenY, LXY1,LXY2,XYY,Yuan} and references therein,
where subsonic flows and nozzles are small perturbations of some
given background flows and nozzles with simple geometries,
respectively.

In this paper, we study the existence of global steady subsonic
Euler flows through periodic nozzles. Consider 2-D steady isentropic Euler equations
\begin{eqnarray}
&&(\rho u)_{x_1}+(\rho v)_{x_2}=0,\label{continuityeq}\\
&&(\rho u^2)_{x_1}+(\rho uv)_{x_2}+p_{x_1}=0,\label{momentumeq1}\\
&&(\rho uv)_{x_1}+(\rho v^2)_{x_2}+p_{x_2}=0,\label{momentumeq2}
\end{eqnarray}
where $\rho$, $(u,v)$, and $p=p(\rho)$ denote the density, velocity,
and  pressure,  respectively. In general, it is assumed that
$p'(\rho)>0$ for $\rho>0$ and $p''(\rho)\geq 0$, where
$c(\rho)=\sqrt{p'(\rho)}$ is called the sound speed. The most
important examples include polytropic gases and isothermal gases.
For polytropic gases, $p=A\rho^{\gamma}$ where $A$ is a constant and
$\gamma$ is the adiabatic constant with $\gamma>1$; and for
isothermal gases, $p=c^2\rho$ with constant sound speed $c$
\cite{CF}.

We consider flows through an infinitely long periodic nozzle given
by
\begin{equation*}
\Omega=\{(x_1,x_2)|f_1(x_1)<x_2<f_2(x_1),-\infty<x_1<\infty\},
\end{equation*}
where $f_{i}$ ($i=1, 2$) is $L-$periodic, i.e., $f_i(x_1+L)=f_i(x_1)$ for $x_1\in \mbR$.
Suppose that
there exist $\alpha\in (0,1)$ and $C>0$ such that
\begin{equation}\label{Eboundary3}
\|f_i\|_{C^{2,\alpha}(\mathbb{R})}\leq C\quad \text{and}\quad \inf_{x_1\in [0, L]}(f_2(x_1)-f_1(x_1))>0.
\end{equation}
Therefore, the domain $\Omega$ satisfies the uniform exterior sphere
condition with some uniform radius $r>0$. Without loss of generality, we assume that $f_1(0)=0$ and $f_2(0)=1$.

Suppose that the nozzle walls are impermeable so that
the flow satisfies the no flow boundary condition
\begin{equation}\label{noflowbc}
(u,v)\cdot \vec{\nu}=0\,\,\text{on}\,\,\partial\Omega,
\end{equation}
where $\vec{\nu}$ is the unit outward normal to the nozzle wall. It
follows from (\ref{continuityeq}) and (\ref{noflowbc}) that
\begin{equation}\label{massflux}
\int_{l}(\rho u,\rho v)\cdot \vec{n} dl \equiv m
\end{equation}
holds for some constant $m$, which is called the mass flux, where
$l$ is any curve transversal to the $x_1-$direction, and $\vec{n}$
is the normal of $l$ in the positive $x_1$-axis direction.

Using the continuity equation \eqref{continuityeq}, when the flow is away from the
vacuum, the momentum equations \eqref{momentumeq1} and \eqref{momentumeq2} are equivalent to
\begin{eqnarray}
uu_{x_1}+vu_{x_2}+h(\rho)_{x_1}=0\label{Enonconservemeq1},\\
uv_{x_1}+vv_{x_2}+h(\rho)_{x_2}=0\label{Enonconservemeq2},
\end{eqnarray}
where $h(\rho)$ is the enthalpy of the flow satisfying
$h'(\rho)=p'(\rho)/\rho$ and can be determined up to a
constant. In this paper, for example, we always choose $h(0)=0$ for
polytropic gases and $h(1)=0$ for isothermal gases. After
determining this integral constant, we denote
$H_0=\inf_{\rho>0}h(\rho)$.

It follows from (\ref{Enonconservemeq1}) and
(\ref{Enonconservemeq2}) that
\begin{equation}\label{EstreamBernoulli}
(u,v)\cdot \nabla (h(\rho)+\frac{1}{2}(u^2+v^2))=0.
\end{equation}
This implies that $\frac{u^2+v^2}{2}+h(\rho)$, which is called
Bernoulli's function, is a constant along each
streamline. For Euler flows in the nozzle, we assume that at  $x_1=0$, Bernoulli function is given,
i.e.,
\begin{equation}\label{Bernoulli}
\left(\frac{u^2+v^2}{2}+h(\rho)\right)(0, x_2)=B_0(x_2),
\end{equation}
where $B_0(x_2)$ is a function defined on $[0,1]$.

When the Bernoulli function $B$ is a  constant, Proposition \ref{Equivalent} shows that the flow is irrotational. The existence of periodic potential flows  with
small mass flux in periodic nozzles was obtained in \cite{LI}. In
this paper, we first study the subsonic and subsonic-sonic
periodic potential flows with relatively large mass flux.
\begin{thm}\label{Thmpotential}
If $B_0(x_2)\equiv \bar{B}>H_0$, then
\begin{enumerate}
\item[1.]
there exists an $\hm>0$, such that for any $m\in(0,\hm)$ there
exists a unique subsonic periodic flow $(\bar\rho, \bar u, \bar v)$ which
satisfies $\inf_{\bar{\Omega}}\bu>0$;
\item[2.]  the maximum of Mach numbers of the flows increases  as $m$ increases and goes to one as $m \to \hm$, i.e., the flows approach sonic;
\item[3.]  there exist a sequence $m_n\rightarrow \hm$ such that the associated potential flows  $(\brho_n, \bu_n, \bv_n)$ converge to $(\hrho, \hu, \hv)$  almost everywhere, which satisfies
\begin{equation}\label{soniceq}
\left\{
\begin{aligned}
\nabla \times (\hu, \hv) =0,\\
\Div(\hrho \hu, \hrho \hv) =0,
\end{aligned}
\right.
\end{equation}
and  the boundary condition \eqref{noflowbc} in the sense of divergence measure field, where $\hrho$ is determined by $\hu$ and $\hv$ via Bernoulli law.
\end{enumerate}
\end{thm}

When $B_0$ is not a constant, we have the following results on subsonic Euler flows in periodic nozzles.
\begin{thm}\label{Maintheorem}
Let the nozzle satisfy (\ref{Eboundary3}) and $B_0$ in
(\ref{Bernoulli}) satisfy
\begin{equation}\label{Bbdycd}
B'_0(0)\geq0,\quad B'_0(1)\leq0.
\end{equation}
For any $m \in (0, \hm)$, there exists $\epsilon_0>0$
such that  if  $B_0(x_2)$
satisfies
\begin{equation}\label{AssumptiononBernoulli}
\|B_0-\bar{B}\|_{C^{1,1}([0,1])}=\epsilon\leq\epsilon_0,
\end{equation}
where $\bar{B}$ is the constant in Theorem \ref{Thmpotential}, then
\begin{enumerate}
\item[1.] (Existence) there exists a periodic flow, i.e.,
$$\rho(x_1+L,x_2)=\rho(x_1,x_2), u(x_1+L,x_2)=u(x_1,x_2), v(x_1+L,x_2)=v(x_1,x_2),$$
which satisfies the original Euler equations
(\ref{continuityeq})-(\ref{momentumeq2}), the boundary condition
(\ref{noflowbc}),   mass flux condition (\ref{massflux}), and the
condition (\ref{Bernoulli});
\item[2.]  (Subsonic flows and positivity of horizontal velocity)  the
flow is globally uniformly subsonic and has positive horizontal
velocity. More precisely, there are two positive constants $c_1$ and $c_2$ such that
\begin{equation}\label{uniformcond}
\sup_{\bar{\Omega}}(u^2+v^2-c^2(\rho))<-c_1 \quad \text{and}\quad
u>c_2;
\end{equation}
\item[3.]  (Regularity)  the flow satisfies
\begin{equation*}
\|\rho\|_{C^{1,\alpha}(\Omega)},\|u\|_{C^{1,\alpha}(\Omega)},
\|v\|_{C^{1,\alpha}(\Omega)}\leq C
\end{equation*}
for some constant $C>0$;
\item[4.] (Uniqueness)
the flow is unique in the class of the periodic flows satisfying
\eqref{uniformcond}.
\end{enumerate}
\end{thm}

A remark about Theorem \ref{Maintheorem} is as follows.

\begin{rmk}
{\rm Using the analysis in this paper, it is easy to show that
there exists a subsonic full compressible Euler flow in the
nozzle, if the entropy is also prescribed at $x_1=0$.}
\end{rmk}

The rest of the paper is organized as follows: In Section
\ref{Reformulation}, we introduce the stream function formulation
for the Euler equations and give the proof of Theorem
\ref{Thmpotential}. In Section \ref{Secstreampb}, a boundary value
problem for stream function is analyzed. This is divided into two
steps. The existence of solutions for the associated problem is
studied in Section \ref{Existence}. The uniqueness, periodicity,
and positivity of horizontal velocity of the flows are proved in
Section \ref{Periodicity}. In Section \ref{SecEuler}, we use the
fixed point theorem to show the existence of Euler flows; the
uniqueness of these flows are obtained by the energy method.


\section{Stream Function Formulation of the Euler Flows and Proof of Theorem \ref{Thmpotential}}\label{Reformulation}
\subsection{Bernoulli's Law and A New Formulation of the Euler Equations.}
We recall that the steady Euler system for subsonic flows is a hyperbolic-elliptic coupled
system \cite{XX3}. Therefore,  one
has to resolve the hyperbolic mode.

To overcome the difficulties mentioned above, we introduce the
stream function for 2-D steady compressible Euler flows so that the
Bernoulli function can be reduced to a single-valued function of
steam function.
This gives an equivalent formulation for Euler flows in terms of
stream function.

\begin{prop}\label{Equivalent}
For a smooth flow away from vacuum with no stagnation point, i.e.,
$u^2+v^2>0$, the Euler system \eqref{continuityeq}-\eqref{momentumeq2} is equivalent to the
 system of equations \eqref{continuityeq}, \eqref{EstreamBernoulli}, and
\begin{eqnarray}\label{Vorticity4}
\omega=\frac{v\partial_{x_1}B-u\partial_{x_2}B}{u^2+v^2},
\end{eqnarray}
where $B=\frac{1}{2}(u^2+v^2)+h(\rho)$ and
$\omega=\partial_{x_1}v-\partial_{x_2}u$ are Bernoulli function and vorticity, respectively.
\end{prop}
\begin{pf}
Let us first show that \eqref{continuityeq}-\eqref{momentumeq2} implies \eqref{Vorticity4}. Differentiating the Bernoulli function with respect to $x_1$ and $x_2$, respectively,  gives
\begin{equation}\label{derBer}
\partial_1 B = u\partial_1 u+v\partial_1 v+\partial_1 h(\rho),\quad  \partial_2 B = u\partial_2 u+v\partial_2 v+\partial_2 h(\rho).
\end{equation}
This, together with (\ref{Enonconservemeq1})-(\ref{Enonconservemeq2}),
yields
\begin{equation}\label{BVorticity2}
\partial_{1}B=v(\partial_1 v-\partial_2 u)=v\omega,\quad
\partial_{2}B=-u(\partial_1 v-\partial_2 u)=-u\omega.
\end{equation}
Therefore,  the equation \eqref{Vorticity4} holds provided $u^2+v^2>0$.

Conversely, it follows from straightforward computation that \eqref{EstreamBernoulli} and \eqref{Vorticity4}  imply \eqref{BVorticity2}. Substituting \eqref{BVorticity2} into \eqref{derBer} gives \eqref{Enonconservemeq1} and \eqref{Enonconservemeq2}. Using \eqref{continuityeq}, one has \eqref{momentumeq1} and \eqref{momentumeq2}.

This finishes the proof of the proposition.
\end{pf}

The continuity equation \eqref{continuityeq} implies that there exists
a stream function $\psi$ such that
$$\partial_{x_1}\psi=-\rho v,\,\ \partial_{x_2}\psi=\rho u.$$
Hence, for the flows away from the vacuum, \eqref{EstreamBernoulli} is
equivalent to
\begin{equation*}
\nabla^{\bot}\psi\cdot\nabla B=0,\label{SB}
\end{equation*}
where $\nabla^{\bot}=(-\partial_{x_2},\partial_{x_1}).$ This
yields that $B$ and $\psi$ are functionally dependent. Therefore,
one may regard $B$ as a function of $\psi$. We denote this
function by $B=\mathcal{B}(\psi)$. It follows from the no flow
boundary condition (\ref{noflowbc}),  that the nozzle walls are
streamlines, so $\psi$ is constant on each nozzle wall. Taking
(\ref{massflux}) into account, one may assume that
\begin{equation}\label{streamBC}
\psi=0\,\,\text{on} \,\,S_1,\,\,\text{and}\,\, \psi=m\,\,\text{on}
\,\,S_2,
\end{equation}
where $S_i =\{(x_1, f_i(x_1))|x_1\in\mbR\}$ ($i=1$, $2$).


As in \cite{XX3}, for any given $s>H_0$, there exist
$\bar{\varrho}=\bar{\varrho}(s)$, $\varrho=\varrho(s)$ and
$\Gamma=\Gamma(s)$ such that
\begin{equation*}\label{Edefmaxcrdensity}
h(\bar{\varrho}(s))=s,\,\,
h(\varrho(s))+\frac{\Gamma^2(s)}{2}=s,\,\,\text{and}\,\,
p'(\varrho(s))=\Gamma^2(s),
\end{equation*}
where $\bar{\varrho}(s)$, $\varrho(s)$, and $\Gamma(s)$ are the
maximum density, the critical density, and the critical speed,
respectively for the states with given Bernoulli constant $s$.
Set
\begin{equation}\label{defSigma}
\Sigma(s)=\varrho(s)\sqrt{2(s-h(\varrho(s)))}.
\end{equation}
Then the straightforward calculations show that
\begin{equation*}
\frac{d\bar{\varrho}}{d s}>0,\quad \frac{d
\varrho}{d s}>0,\quad \text{and}\,\,\,\, \frac{d\Sigma}{d s}>0.
\end{equation*}
Obviously, $\varrho(s)<\bar{\varrho}(s)$, if $s>H_0=\inf_{\rho>0}
h(\rho)$. Furthermore, there exists a $\bar\delta>0$ such that
\begin{equation}\label{rhobound}
\bar{\varrho}(s_1)\geq \varrho(s_2)\quad \text{for any}\,\,\,\,
s_1, s_2\in (\bar B-\bar\delta, \bar B+\bar \delta).
\end{equation}

For a fixed $s$, if $\rho$ and $\mcM$ satisfies
\begin{equation}\label{gBernoulli}
h(\rho) +\frac{\mcM}{2\rho^2} =s,
\end{equation}
then $\rho$ is a two-valued function of $\mathcal{M}$ for $\mathcal{M}\in (0, \Sigma^2(s))$ and the
subsonic branch satisfies $\rho>\varrho(s)$, see \cite{XX3}.  When
$s$ varies,  the subsonic branch will be denoted by
\begin{equation}\label{Edefsubsonicbranch}
\rho=H(\mathcal{M}, s)\,\,\text{for}\,\,(\mathcal{M},s)\in
\{(\mathcal{M},s)|\mathcal{M}\in(0,\Sigma^2(s)), s>H_0\}.
\end{equation}
In view of \eqref{gBernoulli}, we have
\begin{eqnarray*}
\frac{\partial H}{\partial s}=\frac{H^3}{H^2c^2-\mathcal{M}}>0,\quad
\frac{\partial H}{\partial
\mathcal{M}}=\frac{H}{2(\mathcal{M}-H^2c^2)}<0.
\end{eqnarray*}


\subsection{Potential Flows and Proof of Theorem \ref{Thmpotential}.}
If $B_0(x)\equiv \bB$, we have $\mcB(\psi)=\bB$. It follows from \eqref{Edefsubsonicbranch} that $\rho=H(|\nabla\psi|^2, \bB)$. Furthermore, \eqref{Vorticity4} implies $\omega\equiv 0$ because $B(x_1, x_2)\equiv \bB $ in $\Omega$. Therefore, $\psi$ satisfies
\begin{equation}\label{potentialeq}
\Div\left(\frac{\nabla \psi}{H(|\nabla\psi|^2, \bar{B})}\right) =0.
\end{equation}

Let $\{M_n\}$ be a strictly increasing sequence satisfying $\lim_{n\to \infty}M_n =\Sigma^2(\bB)$. We define $H_n (\cdot, \bB)\in C^{\infty} (\mbR)$ satisfying
\begin{equation}
H_n(M, \bB) = \left\{
\begin{aligned}
&H(M, \bB)\quad \text{if}\quad M\leq M_n,\\
&H(M_n, \bB) \quad \text{if}\quad M\geq \frac{M_n+\Sigma^2(\bB)}{2}.
\end{aligned}
\right.
\end{equation}
Combining Lemmas 2.1 and 3.1 in \cite{XX1} and Lemma 1 in \cite{LI} (or Theorem 2 in \cite{Bergner}), there exists a unique periodic solution $\bpsi(\cdot; t)$ of the problem
\begin{equation*}
\left\{
\begin{aligned}
&\Div\left(\frac{\nabla \psi}{H_n(|\nabla\psi|^2, \bar{B})}\right) =0,\\
&\psi=0\quad \text{on}\quad S_1,\quad \psi=t\quad \text{on}\quad S_2
\end{aligned}
\right.
\end{equation*}
for $t\geq 0$. Define $\msM_n(t)=\sup_{\bar{\Omega}}|\nabla
\bpsi_n(\cdot; t)|^2$. Then that $\msM_n(t)$ is a continuous
function of $t$ follows from the same argument in Lemma 4.1 in \cite{XX1}. Let
$m_n=\sup\{t|\msM_n(t)<M_n\}$ and $\hm=\sup_n m_n$. Then as $m\to
\hm$, the maximum of $|\nabla \bpsi|$ of  solutions of the problem
\eqref{potentialeq} and \eqref{streamBC} tend to $\Sigma(\bB)$,
i.e., the flows approach the sonic state. Moreover,  Using the
compensated compactness framework in \cite{CDSW} (Theorem 2.1) and
\cite{XX1} (Theorem 5.1), there exist a sequence $\{m_n\}\uparrow
\bm$ such that the associated $\bu_n=\frac{\partial_2 \bpsi_n}{H(|\nabla\bpsi_n|^2, \bB)}$, $\bv_n = -\frac{\partial_1 \bpsi_n}{H(|\nabla\bpsi_n|^2, \bB)}$, and $\brho_n=H(|\nabla\bpsi|^2, \bB)$ satisfy $\bu_n
\to \hu$, $\bv_n\to \hv$, and $\brho_n\to \hrho$ a.e., where $(\hu, \hv)$ satisfies
$\sup_{\bar{\Omega}}\frac{\hu^2+\hv^2}{p'(\hrho)}=1$, the system \eqref{soniceq},
and the boundary condition \eqref{noflowbc} in the sense of divergence measure field.

The proof of Theorem \ref{Thmpotential} finishes after we prove the following lemma on the properties of subsonic potential flows.
\begin{lem}\label{lemPotentialflow}
If $m>0$, then the solution $\bar{\psi}$ of the problem
\eqref{potentialeq} and \eqref{streamBC} satisfies
$\inf_{\bar{\Omega}}\partial_2\bar{\psi}>0$. Furthermore, as $m$
increases, $\max_{\bar{\Omega}}|\nabla\bpsi|$ also increases.
\end{lem}
\begin{pf}
Note that $\bpsi$ satisfies $0\leq \bpsi\leq m$, therefore,
$\bar{\psi}$ achieves its  minimum and maximum on the boundaries
$S_1$ and $S_2$ respectively. It follows from the Hopf lemma
(Lemma 3.4 in \cite{GT}) that $\partial_{x_2}\bar{\psi}>0$ on
$\partial\Omega$. Therefore, the continuity of
$\partial_{x_2}\bpsi$ yields
\begin{equation*}
\inf_{\substack{(x_1, x_2)\in \partial \Omega,\\ 0\leq x_1 \leq L}}\partial_{x_2}\bar\psi>0.
\end{equation*}
If $\bar\psi$ satisfies \eqref{potentialeq}, then $\partial_{x_2}\bar\psi$ satisfies
\begin{equation}\label{depotentialeq}
\partial_i(\bar{a}_{ij}\partial_j (\partial_{x_2}\bar{\psi}))=0,\quad \text{where}\quad \bar{a}_{ij} =\frac{(H^2c^2-|\nabla \bpsi|^2)\delta_{ij}+\partial_i\bpsi\partial_j\bpsi}{H(|\nabla\bpsi|^2-H^2c^2)}.
\end{equation}
Here and later on, the repeated index means summation from $1$ to
$2$. It follows from the strong maximum principle (Theorem 3.5 in
\cite{GT})  that $\partial_2\bar{\psi}\geq 0$ on $\bar{\Omega}$.
In fact, if the minimum of $\partial_{x_2}\bpsi$ is achieved at
some point $(x_1^*,x_2^*)\in\Omega$, by periodicity, we can always
assume $|x_1^*|\leq L$. Then, it contradicts with the strong
maximum principle for the equation \eqref{depotentialeq} in the
domain $\{(x_1,x_2)\in\Omega:|x_1|\leq \frac{3L}{2}\}$.

It follows from Bernstein estimates (Theorem 15.1 in \cite{GT}) and periodicity of the solutions that
\begin{equation}
\max_{\bar{\Omega}}|\nabla\bpsi|\leq \max_{\partial{\Omega}}|\nabla\bpsi|.
\end{equation}
Thus using comparison principle, Hopf Lemma and periodicity, we can show that the maximum of the flow speed increases as $m$ increases (cf. Lemma 4.4 in \cite{XX1}).
\end{pf}

\subsection{Stream Function Formulation of General Euler Flows.} For any $m\in (0, \hm)$, there exists a unique periodic solution $\bpsi(x_1, x_2)\in C^{2, \alpha}(\bar{\Omega})$ for the problem \eqref{potentialeq} and \eqref{streamBC}. Furthermore, there exist positive constants $\sigma_0$ and $\sigma_1$ such that
\begin{equation}\label{defsigma}
0< \inf_{\bar{\Omega}} \bpsi_{x_2}=\sigma_0\leq \sigma_1
=\sup_{\bar{\Omega}} |\nabla \bpsi|.
\end{equation}
Given $W\in \mcS$ defined by
\begin{equation*}
\mcS=\{ W\in C^{1, \beta}([0, 1]),\,\,\int_0^1 W(s)ds =m,\,\,
\|W-\bpsi_{x_2}(0, \cdot)\|_{C^{1, \beta}[0,1]}\leq \sigma_0/2\},
\end{equation*}
where $\beta\in (0, \alpha)$,  then $W(s)>\sigma_0/2$ for $s\in [0, 1]$. Therefore there exist a function $y=\kappa(\psi)$ such that
\begin{equation}\label{defkappa}
\psi =\int_0^{\kappa(\psi)}  W(s)ds.
\end{equation}
Differentiating \eqref{defkappa} with respect to $\psi$ yields
\begin{equation}\label{derk}
\kappa'(\psi)=\frac{1}{W(\kappa(\psi))}.
\end{equation}
This shows that $\kappa(\psi)\in C^{2, \beta}([0, m])$. Suppose
that the Bernoulli function at $x_1=0$ is $B_0(\kappa(\psi))$ in
terms of $\psi$, since the Bernoulli function is a constant
along each stream line, we get the Bernoulli function defined in
the whole domain $\Omega$ by $B(x_1, x_2) =\mcB(\psi(x_1,
x_2))=B_0(\kappa(\psi(x_1, x_2)))$.  Define
$\mcB(\psi)=B_0(\kappa(\psi))$. Then $\mcB(\psi) \in C^2([0, m])$.
Combining   \eqref{Bbdycd} and \eqref{derk} gives
\begin{equation}\label{endBcond}
\|\mcB-\bar{B}\|_{C^{1,1}([0, m])}\leq C\epsilon, \quad \mcB'(0)\geq
0\quad \text{and}\quad \mcB'(m)\leq 0.
\end{equation}
Therefore,  $\rho=H(|\nabla \psi|^2, \mcB(\psi))$ follows from \eqref{Edefsubsonicbranch} for subsonic flows.  The equation (\ref{Vorticity4}) becomes the following second order equation for the stream function $\psi$,
\begin{equation}\label{Sformulation1}
\Div\left(\frac{\nabla \psi}{H(|\nabla\psi|^2, \mathcal{B}(\psi))}\right)=H(|\nabla \psi|^2, \mcB(\psi))\mathcal{B}'(\psi).
\end{equation}

We first solve the equation
(\ref{Sformulation1}) with the boundary condition  \eqref{streamBC}. Second, we define a map from $W$ to $\psi_{x_2}(0, \cdot)$. The fixed point of this map will induce the existence of the solutions of Euler equations.


\section{Analysis of the Boundary Value Problem for Stream
Function}\label{Secstreampb}

\subsection{Existence of Solutions}\label{Existence}

There are two main difficulties to solve the equation
(\ref{Sformulation1}) in $\Omega$. The first difficulty is that the equation
(\ref{Sformulation1}) becomes degenerate elliptic at sonic
states. In addition, $H$ is not well-defined for arbitrary $\psi$
and $|\nabla\psi|$; neither is $\mathcal{B}$. The second difficulty
is that this is a problem in an unbounded domain. Our basic strategy
is that we extend the definition of $\mathcal{B}$ appropriately,
truncate $|\nabla\psi|$ appeared in $H$ in a suitable way, and use a
sequence of problems in bounded domains to approximate the
problem (\ref{Sformulation1}) and  \eqref{streamBC}.

Set
\begin{equation*} \tilde{g}(s)=\left\{
\begin{array}{ll}
\mathcal{B}'(s),\,\,&\text{if}\,\,0\leq s\leq m,\\
\mathcal{B}'(m)(2m-s)/m,\,\, &\text{if}\,\, m\leq s\leq 2m,\\
\mathcal{B}'(0)(s+m)/m,\,\, &\text{if}\,\, -m\leq s\leq 0,\\
0,\,\,&\text{if}\,\, s\geq 2m\,\,\text{or}\,\, s\leq -m. \\
\end{array}
\right.
\end{equation*}
It is obvious that $\tilde{g}\in C^{0,1}(\mathbb{R})$ and
\begin{equation*}
\|\tilde{g}(s)\|_{C^{0,1}(\mathbb{R}^1)}\leq
\|\mathcal{B}'(s)\|_{C^{0,1}([0,m])}\leq 2\epsilon/\sigma_0.
\end{equation*}
Define
\begin{equation*}
\tilde{\mathcal{B}}(s)=\mathcal{B}(0)+\int_0^s\tilde{g}(t)dt.
\end{equation*}
Then, $\|\tilde{\mathcal{B}}'\|_{C^{0,1}(\mathbb{R}^1)} =
\|\tilde{g}\|_{C^{0,1}(\mathbb{R}^1)}
\leq\|\mathcal{B}'\|_{C^1([0,m])}\leq C\epsilon$. Therefore,
\begin{equation*}
| \tilde{\mcB}(\psi)-\bB|\leq C\epsilon.
\end{equation*}
Hence, there exists $\tilde{\epsilon}_0>0$ such that if
$0<\epsilon<\tilde{\epsilon}_0$, then
$\tilde{\mcB}(\psi)>H_0=\inf_s h(s)$. In view of \eqref{endBcond},
$\tilde{\mcB}$ also satisfies
\begin{equation}\label{tmcBbdycd}
\tilde{\mcB}'(s)\geq 0 \quad \text{for } s\leq 0\quad  \text{and}\quad \tilde{\mcB}'(s)\leq 0\text{ for } s\geq m.
\end{equation}

Let $\ckB=\min_{x\in [0, 1]} B_0(x)$. Choose $\theta_0$  to be a fixed positive constant satisfying
\begin{equation}
0< \theta_0\leq \min\{\Sigma^2(\ckB)/2,\Sigma^2(\bar{B})-\sigma^2_1\},
\end{equation}
where $\sigma_1$ is defined in \eqref{defsigma}.
Let $\zeta\in C^{\infty}(\mbR)$ satisfy
\begin{equation*}
\zeta(s)=\left\{
\begin{array}{ll}
s,\,\, &\text{if}\,\, s<-\theta_0/4,\\
-\theta_0/8,\,\,&\text{if}\,\, s\geq
-\theta_0/8,
\end{array}
\right.
\end{equation*}
and define
$\tilde{\rho}=H(\zeta(|\nabla\psi|^2-\Sigma^2(\tilde{\mathcal{B}}(\psi))+\Sigma^2(\tilde{\mathcal{B}}(\psi),\tilde{\mathcal{B}}(\psi))$.
Instead of the equation (\ref{Sformulation1}), we begin with
investigating the  equation,
\begin{equation}\label{FixedBS2}
\begin{array}{ll}
\partial_1(\frac{\partial_1\psi}{\tilde{\rho}})+\partial_2(\frac{\partial_2\psi}{\tilde{\rho}})=\tilde{\rho}\tilde{\mathcal{B}}'(\psi).
\end{array}
\end{equation}
The equation \eqref{FixedBS2} can also be written in the following non-divergence form
\begin{equation}\label{FixedBS3}
A_{ij}(\nabla\psi,\psi)\partial_{ij}\psi=\mathcal{F}(\nabla\psi,\psi)\end{equation}
where
\begin{equation*}
A_{ij}(\nabla\psi,\psi))=\trho\delta_{ij}+\frac{\zeta'\trho}{(\trho^2\tc^2 -(\zeta+\Sigma^2))}\partial_i\psi\partial_j\psi,
\end{equation*}
and
\begin{equation*}
\mathcal{F}(\nabla\psi,
\psi)=\tilde{\mcB}'(\psi)\left(\trho^3
\frac{\trho^2 \tc^2-(\zeta+\Sigma^2)+|\nabla\psi|^2}{\trho^2\tc^2 -(\zeta +\Sigma^2)} + \frac{(\zeta'-1)\Sigma\Sigma'|\nabla \psi|^2}{\trho^2\tc^2 -(\zeta +\Sigma^2)}
\right),
\end{equation*}
where the variables in $\zeta$ and $\Sigma$ are
$|\nabla\psi|^2-\tilde{\mcB}(\psi)$ and $\tilde{\mcB}(\psi)$,
respectively. The direct calculation shows that the eigenvalues
$\Lambda$ and $\lambda$ of $[A_{ij}]_{2\times2}$ satisfy
$C^{-1}\leq |\Lambda/\lambda|\leq C$ for some constant $C$ and
thus the equation  \eqref{FixedBS3} is uniformly elliptic.
However, $\mathcal{F}$ involves a quadratic growth in
$|\nabla\psi|$, so it is not easy to apply the classical elliptic
theory directly. The strategy is to modify $\mcF$ by
\begin{equation}
\tilde{\mcF}(\nabla\psi,
\psi)=\tilde{\mcB}'(\psi)\frac{\trho^5\tc^2}
{\trho^2\tc^2 -(\zeta +\Sigma^2)},
\end{equation}
It is easy to see that $\tilde{\mcF}=\mcF$ when $\psi$ satisfies $|\nabla\psi|^2-\Sigma^2(\tilde{\mcB}(\psi))\leq -\theta_0/4$.

\begin{prop}\label{Fixedexistence}
There exists a solution $\psi(\cdot; t)\in C^{2, \alpha}(\bar{\Omega})$ of the problem
\begin{equation}\label{41}
\left\{
\begin{aligned}
&A_{ij}\partial_{ij}\psi=\tilde{\mcF}\\
&\psi=0 \quad \text{on } S_1, \quad \psi=t\text{ on } S_2.
\end{aligned}
\right.
\end{equation}
Furthermore, if $t\leq m$, then
\begin{equation}\label{firstmaximum}
0\leq\psi\leq m.
\end{equation}
There exists $m_1>0$ and $\epsilon_1>0$ such that if $0\leq t\leq
m_1$ and $\|\mcB'\|\leq \epsilon_1$, then
\begin{equation}\label{Property}
|\nabla\psi|^2-\Sigma(\mcB(\psi)) \leq -\theta_0/2.
\end{equation}
\end{prop}
\begin{pf}
For the unbounded domain $\Omega$,  we use a sequence of
boundary value problem defined in bounded domains $\tOmegaN$ to
approximate the problem in $\Omega$, where $\tOmegaN \in C^{2, \alpha}$ satisfies $\Omega_{N}\subset \tOmegaN\subset \Omega_{2N}$ with $\Omega_k :=\Omega\cap\{|x_1|\leq k\}$. The construction of $\tOmegaN$ can be found in \cite{XX1}.

We first solve the boundary value problem
\begin{equation}\label{Modifiedbvp2}
\left\{
\begin{array}{ll}
A_{ij}(\nabla\psi,\psi)\partial_{ij}\psi=\tilde{\mathcal{F}}(\nabla\psi,\psi)\,\
\text{in}\,\ \tOmegaN,\\
\psi=\frac{x_2-f_1(x_1)}{f_2(x_1)-f_1(x_1)}t\,\ \text{on}\,\
\partial\tOmegaN.
\end{array}
\right.
\end{equation}
Similar to the proof of Proposition 3 in \cite{XX3}, there exists a
$C^{2,\alpha}$-solution $\psi_N$, such that
$$|\psi_N|\leq
C\left(t+\left|\frac{\mathcal{\tilde{F}}}{\lambda}\right|_0\right),\,\
\|\psi_N\|_{2,\alpha;\tOmegaN}\leq C\left(
\Lambda/\lambda,|f_i|_{C^{2,\alpha}},m,
\left|\frac{\mathcal{\tilde{F}}}{\lambda}\right|_0\right),
$$
By Arzela-Ascoli theorem, one can select a subsequence of
$\{\psi_N\}$ (We still label it by $\{\psi_N\}$), such that
$$\psi_N\rightarrow\psi\,\ \text{in}\,\ C^{2,\alpha_1}(K)\,\ \text{with}\,\ K\Subset\bar{\Omega},\,\ \alpha_1<\alpha.$$
Obviously, $\psi$  solves (\ref{41}).

Since $\tilde{\mcB}$ satisfies \eqref{tmcBbdycd}, then
$\widetilde{\mcF}(\nabla\psi_N,\psi_N)\geq0$ in the domain
$\tOmegaN\cap\{\psi_N\geq m\}$. Thus according to system
(\ref{Modifiedbvp2}),
$$A_{ij}(\nabla\psi_N,\psi_N)\partial_{ij}\psi_N\geq0\,\ \text{in}\,\ \tOmegaN\cap\{\psi_N\geq m\}.$$
Since
$$\psi_N\leq m\,\ \text{on}\,\ \partial\left(\tOmegaN\cap\{\psi_N\geq
m\}\right),$$ by the maximum principle (Theorem 3.1 in \cite{GT}),
one has
$$\psi_N\leq\sup_{\partial\tOmegaN}\psi_N\leq m\,\ \text{in}\,\ \tOmegaN\cap\{\psi_N\geq m\}.$$
Similarly, it is also true that
$$\psi_N\geq\inf_{\partial\tOmegaN}\psi_N\geq 0\,\ \text{in}\,\ \tOmegaN\cap\{\psi_N\leq 0\}.$$
Therefore, one has $0\leq \psi_N\leq m$. So the limit $\psi$ satisfies \eqref{firstmaximum}.

The H\"{o}lder estimate for the gradients of elliptic equations of two variables implies that
$$[\psi_N]_{1,\mu;\tOmegaN}=\sup_{x,y\in\tOmegaN}
\frac{|\nabla\psi_N(x)-\nabla\psi_N(y)|} {|x-y|^{\mu}} \leq
C(\Lambda/\lambda,|f_i|_2)
\left(1+t+\left|\frac{\tilde{\mathcal{F}}}{\lambda}\right|_0\right).$$
Then, $\psi_N$ satisfies the estimate
$$|\nabla\psi_N|^2\leq \eta\left(1+t+\left|\frac{\tilde{\mathcal{F}}}{\lambda}\right|_0\right)^2+C_{\eta}\left(t+\left|\frac{\tilde{\mathcal{F}}}{\lambda}\right|_0\right)^2,$$
where $C_{\eta}$ is independent of $N$. Note that
\begin{equation}
|\tilde{\mcF}|\leq C\epsilon,
\end{equation}
there exist $\eta_1$, $m_1$ and $\epsilon_1$ such that
\begin{equation}
\eta_1\left(1+m_1+C\epsilon_1\right)^2\leq \min \{\frac{\Sigma^2(\ckB)}{4}, \frac{\sigma_1^2}{2}\}, \quad
C_{\eta_1}\left(m_1+C\epsilon_1\right)^2 \leq \min \{\frac{\Sigma^2(\ckB)}{4}, \frac{\sigma_1^2}{2}\}.
\end{equation}
If $t\leq m_1$, $\epsilon\leq \epsilon_1$, then
\[
|\nabla\psi_N|^2 -\Sigma^2(\mcB(\psi_N))\leq
-\frac{\Sigma^2(\ckB)}{2}\leq -\theta_0.
\]
Thus, the limit $\psi$ satisfies \eqref{Property}  if  $\epsilon_1$ is suitably small.

This finishes the proof of the proposition.
\end{pf}

\subsection{Uniqueness, Periodicity and Positivity of Horizontal Velocity}\label{Periodicity}
In this subsection, we show that the uniformly subsonic solution obtained in Proposition \ref{Fixedexistence} for the problem  (\ref{Sformulation1}) with the boundary condition
\begin{equation}\label{BCt}
\psi=0 \quad \text{on}\quad S_1,\quad \text{and}\quad \psi=t\quad \text{on}\quad S_2.
\end{equation}
is unique.

\begin{prop}\label{uniqueness}
There exists a positive constant $\epsilon_2<\epsilon_1$ such that if $\|\mathcal{B}'\|_{C^{0,1}}=\epsilon\leq\epsilon_2$, then for $0\leq t\leq m$,
the solution $\psi$ of (\ref{Sformulation1}) and \eqref{BCt}
which satisfies
\begin{equation}\label{ellipcond}
|\nabla\psi|^2 -\Sigma(\mcB(\psi)) \leq -\theta_0/4, \quad 0\leq \psi\leq m
\end{equation}
must be unique.
\end{prop}

\begin{pf}
Let $\psi_i, i=1,2$, both solve (\ref{Sformulation1}) and \eqref{BCt}. Set $\Psi=\psi_1-\psi_2$. Then $\Psi$ satisfies
the following elliptic boundary value problem,
\begin{equation}
\left\{
\begin{array}{ll}
\partial_i(a_{ij}\partial_j\Psi)+\partial_i(b_i\Psi)=c_i\partial_i\Psi+d\Psi,\,\
\text{in}\,\ \Omega,\label{Uniqueness1}\\
\Psi=0,\,\ \text{on}\,\ \partial\Omega,
\end{array}
\right.
\end{equation}
where
\begin{eqnarray*}
&&a_{ij}=\int^1_0\frac{(\tH^2\tc^2-|\nabla\tilde{\Psi}|^2)\delta_{ij}+\partial_i\tilde{\Psi}\partial_j\widetilde{\Psi}}{(|\nabla\widetilde{\Psi}|^2-\tH^2\tc^2)\tH}d\theta, \quad b_i=-\int^1_0\frac{\partial_i\widetilde{\Psi}\tH\mathcal{B}'(\widetilde{\Psi})}{\tH^2\tc^2-|\nabla\widetilde{\Psi}|^2}d\theta,\\
&& c_i=\int^1_0\frac{\tH\partial_i
\tilde{\Psi}}{|\nabla\tilde{\Psi}|^2-\tH^2\tc^2}\mcB'(\tilde{\Psi})d\theta,\quad
d=\int^1_0\frac{\tH^3}{\tH^2\tc^2-|\nabla\widetilde{\Psi}|^2}[\mathcal{B}'(\widetilde{\Psi})]^2+\tH\mathcal{B}''(\widetilde{\Psi})d\theta,
\end{eqnarray*}
with $\tc=\sqrt{p'(\tH)}$,
$\tH=H(|\nabla\widetilde{\Psi}|^2,\mathcal{B}(\widetilde{\Psi}))$,
and  $\widetilde{\Psi}=\theta \psi_1+(1-\theta)\psi_2$. Since
$\|\mcB'\|_{C^{0,1}}\leq C\epsilon$,  we have
\begin{equation*}
|b_i|+|c_i| +|d| \leq C\epsilon.
\end{equation*}

Choose a smooth cut-off function $\eta\in C^{\infty}_0(\mbR)$ satisfying
\begin{equation}\label{defeta}
\eta=\eta(s)=\left\{
\begin{array}{ll}
0,\,\ |s|\geq(N+1)L,\\
1,\,\ |s|\leq NL.
\end{array}
\right.
\end{equation}
Multiplying $\eta^2(x_1)\Psi$ on both sides of the equation
(\ref{Uniqueness1}) and integrating by parts yields that
\begin{equation*}
\begin{aligned}
\iint_{\{|x_1|\leq
NL\}\cap \Omega}|\nabla\Psi|^2 dx_1 dx_2&\leq C\iint_{\Omega}\eta^2(a_{ij}\partial_i\Psi\partial_j\Psi) dx_1 dx_2\\
&\leq C\iint_{\Omega}(|\eta\eta'\Psi|(|\Psi|+|\nabla\Psi|)+\epsilon\eta^2|\Psi|(|\Psi|+|\nabla\Psi|)) dx_1dx_2\\
& \leq C\iint_{\{NL\leq|x_1|\leq(N+1)L\}\cap\Omega}(|\Psi|^2+|\nabla\Psi|^2) dx_1dx_2\\
& \quad  +
C\epsilon\iint_{\{|x_1|\leq NL\}\cap \Omega}(|\Psi|^2 +|\nabla \Psi|^2) dx_1 dx_2.
\end{aligned}
\end{equation*}
Since $\Psi=0$ on $\partial\Omega$,  by Poincare's inequality, we have
$$\iint_{\{kL\leq|x_1|\leq(k+1)L\}\cap \Omega}|\Psi|^2 dx_1dx_2\leq C\iint_{\{kL\leq|x_1|\leq(k+1)L\}\cap\Omega}|\nabla\Psi|^2 dx_1 dx_2$$
for any $k\in \mbZ$. When $\epsilon\leq\epsilon_2$ for some
$\epsilon_2>0$, we have
\begin{equation}
\iint_{\{|x_1|\leq NL\}\cap \Omega}|\nabla\Psi|^2 dx_1 dx_2\leq C\iint_{\{NL\leq|x_1|\leq
(N+1)L\}\cap \Omega}|\nabla\Psi|^2 dx_1dx_2. \label{Poincareinequality},
\end{equation}
Noting that $|\nabla\Psi|\in L^{\infty}(\Omega)$, it  yields that
$\nabla\Psi\in L^2(\Omega)$. Thus
$$\lim_{N\rightarrow\infty}\iint_{\{NL\leq|x_1|\leq(N+1)L\}\cap\Omega}|\nabla\Psi|^2dx_1dx_2=0.$$
Therefore, the estimate (\ref{Poincareinequality}) implies that
$\iint_{\Omega}|\nabla\Psi|^2\equiv0$. It follows from $\Psi\equiv0$ on
$\partial\Omega$ that  $\Psi\equiv0$ in $\Omega$. Therefore, the solution of (\ref{Sformulation1}) and \eqref{BCt} is unique.
\end{pf}

One can check easily that if $\psi(x_1,x_2)$ solves the boundary
value problem (\ref{Sformulation1}) and \eqref{streamBC}, so does $\psi(x_1+L,x_2)$. Then
the uniqueness implies the following corollary.

\begin{cor}\label{corper}
For any $\theta_0>0$, there exists a positive constant
$\epsilon_2<\epsilon_1$ such that if
$\|\mathcal{B}'\|_{C^{0,1}}=\epsilon\leq\epsilon_2$, then for $0\leq t\leq m$,  the
solution $\psi$ of (\ref{Sformulation1}) and \eqref{BCt}
satisfying  \eqref{ellipcond} must be periodic with respect to
$x_1$ with period $L$, i.e.,
$$\psi(x_1,x_2)=\psi(x_1+L,x_2), \forall (x_1,x_2)\in\Omega.$$
\end{cor}

Now we show that $\psi_{x_2}(0, \cdot)\in \mcS$.
\begin{prop}\label{Lowerbound}
For any $\theta_0>0$, there exists a positive constant
$\epsilon_3<\epsilon_2$ such that if
$\|\mathcal{B}'\|_{C^{0,1}}=\epsilon\leq\epsilon_3$, then  for $0\leq t\leq m$, the
solution $\psi$ of (\ref{Sformulation1}) and \eqref{BCt}
satisfying \eqref{ellipcond}, then $\|\nabla \psi-\nabla
\bpsi\|_{C^{1, \alpha}}\leq C\epsilon$.
\end{prop}
\begin{pf}
Set $\Psi=\psi-\bar{\psi}$. Subtracting (\ref{potentialeq}) from
(\ref{Sformulation1}) gives that $\Psi$ satisfies
\begin{equation}\label{Lowerbound2}
\left\{
\begin{aligned}
&\partial_i(a_{ij}\partial_j\Psi +b_i (\mcB(\psi)-\bB)) = H(|\nabla\psi|^2, \mcB(\psi))\mcB'(\psi),\quad \text{in}\quad \Omega,\\
&\Psi=0\quad \text{on}\quad \Omega,
\end{aligned}
\right.
\end{equation}
where
\begin{eqnarray*}
a_{ij}=\int^1_0\frac{(\hH^2\hc^2-|\nabla\hPsi|^2)\delta_{ij}+\partial_i\hPsi\partial_j\hPsi}{\hH(\hH^2\hc^2-|\nabla\hPsi|^2)}d\theta,\quad
b_i=-\int^1_0\frac{\partial_i\hPsi\hH}{\hH^2
\hc^2-|\nabla\hPsi|^2}d\theta,
\end{eqnarray*}
with $\hc=\sqrt{p'(\hH)}$,  $\hH=H(|\nabla\hPsi|^2, \theta
\mcB(\psi)+(1-\theta)\bB)$, $\hPsi=\theta\psi+(1-\theta)\bpsi$.
Since both $\psi$ and $\bar{\psi}$ are periodic, $\Psi$ is also
periodic with period $L$.  Multiplying the equation in
\eqref{Lowerbound2} with $\Psi$ and integrating the resulting
equation on $\Omega\cap \{-2L\leq x_1\leq 2L\}$, and  integration
by parts yield
\begin{equation*}
\|\nabla \Psi\|_{L^2(\Omega\cap \{-2L\leq x_1\leq 2L\})}\leq
C\epsilon.
\end{equation*}
Applying Moser's iteration (Theorems 8.17 and 8.25 in \cite{GT}), we have
\begin{equation*}
\|\Psi\|_{L^\infty(\Omega\cap \{-\frac{L}{2}\leq x_1\leq
\frac{3L}{2}\})} \leq C \epsilon.
\end{equation*}
Using the estimate for elliptic equation of two variables  (Theorem 12.4 and global estimates on page 304 in \cite{GT}), we have
\begin{equation*}
\|\Psi\|_{C^{1, \alpha}(\bar\Omega\cap\{0\leq x_1\leq L\})}\leq C
\epsilon.
\end{equation*}
Using Schauder estimate, we have $\|\Psi\|_{C^{2, \alpha}(\Omega_1)}\leq C\epsilon$.
Choosing $\epsilon_3>0$ sufficiently small shows that
$\|\nabla\Psi\|_{C^{1, \alpha}(\Omega_1)}\leq \sigma_0/2$ provided $0< \epsilon\leq \epsilon_3$. This finishes the proof of the
proposition.
\end{pf}

\section{Existence and Uniqueness of the Euler Flows}\label{SecEuler}

In this section,  We  prove the existence and uniqueness of subsonic solutions for Euler equations.

\begin{prop}\label{propexistence}
There exists a positive constant $\epsilon_5$ such
that if $\|B_0'\|_{C^1}\leq \epsilon_5$, then the system
\eqref{continuityeq}-\eqref{momentumeq2} under the conditions
(\ref{noflowbc}),   mass flux condition (\ref{massflux}), and the
condition (\ref{Bernoulli}) has a subsonic solution with positive
horizontal velocity.
\end{prop}
\begin{pf}
It follows from Proposition \ref{Fixedexistence} that there exists a solution $\psi(\cdot; t)$ for the problem  \eqref{41}. Set
\begin{equation*}
S(t)=\{\psi(\cdot; t) |\psi(\cdot;t) \text{ solves } \eqref{41}\}
\end{equation*}
and
\begin{equation*}
\Delta(t)=\inf_{\psi\in S(t)}\sup_{\bar{\Omega}}\{|\nabla \psi(\cdot; t)|^2 -\Sigma^2(\tB(\psi(\cdot; t)))\}
\end{equation*}
Let $\tm =\sup_s\{s\in(0,\hat{m})|\Delta(s)\leq
-\theta_0/4\}$. Then $\Delta(t)$ is continuous  for $t\in [0, \tilde{m}]$(cf. Proposition 6 in
\cite{XX3}). We claim that $\tm > m$. Indeed, it follows from
Proposition \ref{Fixedexistence} that $\Delta(m_1)\leq -\theta_0/2$.
If $\tilde{m}\leq m$, then
\begin{equation*}
\begin{aligned}
&|\nabla \psi(\cdot; \tm)|^2-\Sigma^2(\mcB(\psi(\cdot; \tm)))\\
 =& |\nabla \psi(\cdot; \tm)|^2- |\nabla \bpsi(\cdot; \tm)|^2 + |\nabla \bpsi(\cdot; \tm)|^2-\Sigma^2(\bB)+ \Sigma^2(\bB) - \Sigma^2(\mcB(\psi(\cdot; \tm))) \\
\leq &  C\epsilon -\theta_0 +C\epsilon,
\end{aligned}
\end{equation*}
where we use the property that the maximum of flow speed increases
as the mass flux increases in Lemma \ref{lemPotentialflow}. Thus
there exists a positive $\epsilon_4\leq \epsilon_3$ such that
$|\nabla \psi(\cdot; \tm)|^2-\Sigma^2(\mcB(\psi(\cdot; \tm)))\leq
-\theta_0/2$. This contradicts with $\Delta(\tilde{m})
=-\theta_0/4$ which follows from the continuity of $\Delta(s)$.
Thus  $\tilde{m}>m$. Thus it implies that for $t=m$ the problem
\eqref{41} has a unique
solution satisfying
\begin{equation*}
|\nabla\psi|^2 -\Sigma(\mcB(\psi))\leq -\theta_0/4.
\end{equation*}
Using Proposition \ref{Lowerbound}, there exists a positive
constant $\epsilon_4<\epsilon_3$ such that $|\nabla \psi-\nabla
\bar{\psi}|\leq \sigma_0/2$. In particular,
$|\partial_{x_2}\psi(0, \cdot)-\partial_{x_2}\bpsi(0, \cdot)|\leq
\sigma_0/2$. Therefore $\partial_{x_2}\psi(0, \cdot) \in \mcS$.
 Hence, we can define a map $T: \mcS \rightarrow \mcS$ by
\begin{equation*}
T(W) =\psi_{x_2}(0, \cdot).
\end{equation*}

By Proposition  \ref{Lowerbound}, $\|\psi_{x_2}(0, \cdot)-\bpsi_{x_2}(0,\cdot)
\|_{C^{1, \alpha}}\leq C \epsilon$. Thus $T\mcS$ is a compact
subset of $\mcS$. It is easy to see that $T$ is a continuous map. Note that $\mcS$ is a closed convex set in
$C^{1, \beta}[0, 1]$;  the existence of fixed point of $T$ on
$\mcS$ follows from Schauder fixed point theorem (Theorem 11.1 in
\cite{GT}).

Let $\hat{W}$ be the fixed point of $T$ and $\hat{\kappa}(y)$ satisfies
\begin{equation*}
y =\int_0^{\hat{\kappa}(y)}\hat{W}(s)ds.
\end{equation*}
Define $\hat{\mcB}(\psi)=B_0(\hat{\kappa}(\psi))$. Then
\begin{equation*}
\left\{
\begin{aligned}
& \Div\left(\frac{\nabla \psi}{H(|\nabla\psi|^2, \hat{\mcB}(\psi))}\right)=H(|\nabla\psi|^2, \hat{\mcB}(\psi)) \hat{\mcB}'(\psi)\\
& \psi=0\text{ on } S_1 \text{   and   } \psi=m \text{ on } S_2
\end{aligned}
\right.
\end{equation*}
has a solution satisfies $\hat{\psi}_{x_2}(0, x_2) =\hat{W}(x_2)$. It is clear that $\rho=H(|\nabla \hat{\psi}|^2, \hat{\mcB}(\psi))$, $u =\partial_{x_2}\hat{\psi}/\rho$, and $v=-\partial_{x_2}\hat{\psi}/\rho$ satisfy the Euler equations and the condition
\begin{equation*}
\left(\frac{u^2+v^2}{2} +h(\rho)\right)(0, x_2)=\hat{\mcB}(\psi(0, x_2))=B_0(x_2).
\end{equation*}
Thus we get the existence of solution for Euler system under the conditions
(\ref{noflowbc}),   mass flux condition (\ref{massflux}), and the
condition (\ref{Bernoulli}). Furthermore, it follows from Proposition \ref{Lowerbound}
that horizontal velocity is positive.
\end{pf}

Now we show that the periodic Euler flow is also unique.

\begin{prop}\label{propunique}
There exists a positive constant $\epsilon_6$ such
that if $\|B_0'\|_{C^1}\leq \epsilon_6$, then the uniformly subsonic solution
the system \eqref{continuityeq}-\eqref{momentumeq2} under the
conditions (\ref{noflowbc}),   mass flux condition
(\ref{massflux}), and the condition (\ref{Bernoulli}) is unique.
\end{prop}
\begin{pf}
Suppose that $\psi_i$ ($i=1, 2$) are stream functions of two periodic solutions of the Euler equations with positive velocity. Let $\kappa_i$ ($i=1, 2$) be the functions satisfying  $\kappa_i(\psi_i(0, x_2)) =x_2$ ($i=1, 2$). Set $\mcB_i(\psi)=B_0(\kappa_i(\psi))$. Then $\psi_i$ satisfies the problem
\begin{equation*}
\left\{
\begin{aligned}
&\Div\left(\frac{\nabla \psi}{H(|\nabla\psi|^2, \mcB_i(\psi))}\right) = H(|\nabla\psi|^2, \mcB_i(\psi))\mcB_i'(\psi),\quad \text{in}\quad \Omega,\\
& \psi=0\quad \text{on}\quad S_1, \quad \psi=m\quad \text{on}\quad S_2.
\end{aligned}
\right.
\end{equation*}

Set
$\Psi=\psi_1-\psi_2$. Then $\Psi$ satisfies the problem
\begin{equation}\label{unique2}
\left\{
\begin{aligned}
&\partial_i(a_{ij}\partial_j\Psi)+\partial_i(b_i\mcD)=c_i\partial_i\Psi+d\mcD+e\mcE, \quad \text{in}\quad \Omega,\\
&\Psi=0 \quad \text{on}\quad \partial \Omega,
\end{aligned}
\right.
\end{equation}
where
\begin{eqnarray*}
&&a_{ij}=\int^1_0\frac{(H^2c^2-|\nabla\tilde{\Psi}|^2)\delta_{ij}+\partial_i\tilde{\Psi}\partial_j\tilde{\Psi}}{H(H^2c^2-|\nabla\tilde{\Psi}|^2)}d\theta,
\quad
b_i=-\int^1_0\frac{\partial_i\tilde{\Psi}H}{H^2c^2-|\nabla\tilde{\Psi}|^2}d\theta\\
&&
c_i=\int_0^1 \frac{H\partial_i\tilde{\Psi}}{|\nabla\tilde{\Psi}|^2-H^2c^2}(\theta \mcB_1'(\psi_1)+(1-\theta)\mcB_2'(\psi_2))\\
&&d=\int_0^1 \frac{H^3}{H^2c^2 -|\nabla\tilde{\Psi}|^2}(\theta
\mcB_1'(\psi_1) +(1-\theta)\mcB_2'(\psi_2) d\theta, \quad e=
\int_0^1 Hd\theta,
\end{eqnarray*}
with $\tilde{\Psi}=\theta\psi_1+(1-\theta)\psi_2,
H=H(|\nabla\tilde{\Psi}|^2, \theta
\mcB_1(\psi_1)+(1-\theta)\mcB_2(\psi_2)$,
$\mcD=\mathcal{B}_1(\psi_1)-\mathcal{B}_2(\psi_2)$ and
$\mcE=\mathcal{B}'_1(\psi_1)-\mathcal{B}'_2(\psi_2)$.

We first multiply $\Psi$ on both sides of the equation in
\eqref{unique2} and integrate the resulting equation  on
$\Omega_1$. In view of periodicity of the coefficients and $\Psi$
in \eqref{unique2}, integration by parts yields
\begin{equation}
\int_{\Omega_1} a_{ij} \partial_i\Psi\partial_{j}\Psi dx=\int_{\Omega_1}(b_i\mcD \partial_i\Psi+ c_i\Psi \partial_i\Psi) +d \mcD\Psi
+e\mcE\Psi dx.
\end{equation}
Thus,
\begin{equation}\label{unique3}
\begin{array}{lll}
\int_{\Omega_1}|\nabla\Psi|^2 dx &\leq&
C\epsilon \int_{\Omega_1}(|\Psi|^2+|\nabla\Psi|^2)+C\int_{\Omega_1}(|\mcD|^2+|\mcE|^2)dx.
\end{array}
\end{equation}
Since $\Psi\in W^{1,\infty}$, the first term on the right hand
side of (\ref{unique3}) is uniformly bounded.

Note that $\Psi=0$ on $\partial\Omega\cap \bar{\Omega_1}$,
Poincare inequality implies that
$$\int_{\Omega_1}|\Psi|^2 dx \leq C\int_{\Omega_1}|\nabla\Psi|^2dx.$$ Therefore, the first term on the right hand side of \eqref{unique3} can be
absorbed by the left hand side.

Now let us estimate the second term on the right hand of \eqref{unique3}. First,
\begin{equation*}
\begin{aligned}
\int_{\Omega_1}|\mcD|^2 dx  = &\int_{\Omega_1}|B_0\circ \kappa_{1}(\psi_1)-B_0\circ
\kappa_{2}(\psi_2)|^2 dx\\
\leq &\int_{\Omega_1}|B_0\circ \kappa_{1}(\psi_1)-B_0\circ
\kappa_{1}(\psi_2)|^2 dx+\int_{\Omega_1}|B_0\circ
\kappa_{1}(\psi_2)-B_0\circ \kappa_{2}(\psi_2)|^2 dx\\
=& I_1+I_2.
\end{aligned}
\end{equation*}
Using mean value theorem, we have
$$I_1\leq C\epsilon \int_{\Omega_1}|\Psi|^2 dx \leq C \epsilon \int_{\Omega_1}|\nabla\Psi|^2 dx.$$
Note that $\psi_1(0, \kappa_1(\psi_2)) =\psi_2(0, \kappa_2(\psi_2))$, we have
\begin{equation*}
\int_0^{\kappa_1(\psi_2)}\partial_{x_2}\psi_1(0, s)-\partial_{x_2}\psi_2(0, s)ds =\int_{\kappa_1(\psi_2)}^{\kappa_2(\psi_2)}\partial_{x_2}\psi_2(0, s)ds.
\end{equation*}
In view of the fact that $\partial_{x_2}\psi_2\geq
\frac{\sigma_0}{2}$, we have
\begin{equation*}
|\kappa_1(\psi_2)-\kappa_2(\psi_2)|\leq C\|\nabla\Psi(0, \cdot)\|_{L^2[0,1]}
\end{equation*}
Thus,
\begin{equation*}
\begin{array}{lllll}
I_2&\leq& \int_{\Omega_1} C \epsilon \|\nabla \Psi(0,\cdot)\|_{L^2[0,1]}dx_1dx_2\leq C \epsilon \|\nabla \Psi\|_{L^{\infty}(\Omega_1)}.
\end{array}
\end{equation*}
Similarly, we can show that $\int_{\Omega_1}|\mcE|^2dx_1dx_2\leq C\epsilon\|\nabla\Psi\|_{L^{\infty}(\Omega_1)}$.

This implies that
\begin{equation*}
\|\nabla \Psi\|_{L^2(\Omega_1)}\leq C\epsilon \|\nabla
\Psi\|_{L^{\infty}(\Omega_1)}.
\end{equation*}
By Nash-Moser's iteration, we have
\begin{equation*}
\|\Psi\|_{L^{\infty}(\Omega_1)}\leq C\epsilon \|\nabla \Psi\|_{L^{\infty}(\Omega_1)}
\end{equation*}
Using the estimate for elliptic equation of two variables (Theorem 12.4 and global estimates on page 304 in \cite{GT}), we have
\begin{equation*}
\|\Psi\|_{C^{1, \alpha}(\Omega_1)}\leq C \epsilon
\|\Psi\|_{L^{\infty}(\Omega_1)}= C \epsilon
\|\Psi\|_{L^{\infty}(\Omega_1)}
\end{equation*}
Therefore $\Psi\equiv 0$ in $\Omega$.

This finishes the proof of the proposition.
\end{pf}

Choosing $\epsilon_0 =\min\{\epsilon_1, \epsilon_2, \epsilon_3,
\epsilon_5, \epsilon_6\}$, then Theorem \ref{Maintheorem} follows
from Propositions \ref{propexistence} and \ref{propunique}.

\bigskip
{\bf Acknowledgement.} The second author thanks Professor Louis
Nirenberg for proposing this problem and  helpful discussions.
Part of the work was done when the second author was visiting The
Institute of Mathematical Sciences, The Chinese University of Hong
Kong. He thanks the institute for its hospitality and support. Both
of the authors thank Professor Zhouping Xin for helpful
discussions.

\end{document}